\documentclass[12pt,psamsfonts]{amsart}

\RequirePackage{amsfonts}

 \RequirePackage{amssymb}
\RequirePackage{latexsym}

\begin{document}
\thispagestyle{empty} \setcounter{page}{1}

\title[Enhanced delay to bifurcation]
{Enhanced delay to bifurcation}

\author[Jean--Pierre Fran\c{c}oise, Claude Piquet  and Alexandre Vidal]
{Jean--Pierre Fran\c{c}oise, Claude Piquet  and Alexandre Vidal}

\address{Laboratoire J.-L. Lions, UMR 7598, CNRS, Universit\'e P.-M. Curie, Paris6,
Paris, France}

\email{Jean-Pierre.Francoise@upmc.fr}

\thanks{This research was supported by the ACI ``Instantbio"}


\date{}
\dedicatory{}

\maketitle

\maketitle

\begin{abstract}
We present an example of slow-fast system which
displays a full open set of initial data so that the corresponding
orbit has the property that given any $\epsilon$ and $T$, it
remains to a distance less than $\epsilon$ from a repulsive part
of the fast dynamics and for a time larger than $T$. This example
shows that the common representation of generic fast-slow systems
where general orbits are pieces of slow motions near the
attractive parts of the critical manifold intertwined by fast
motions is false. Such a description is indeed based on the
condition that the singularities of the critical set are folds. In our
example, these singularities are transcritical.
\end{abstract}

\section{Introduction}

A first approximation for the time evolution of fast-slow systems
dynamics is often seen as follows. A generic orbit jumps close to
an attractive part of the equilibrium set of the fast dynamics. It
evolves then slowly close to this attractive part until that,
 under the influence of the slow dynamics, this attractive part
bifurcates into a repulsive one. Then the generic orbit jumps to
another attractive part of the fast equilibrium set until it
 also looses its stability and then jumps again in search of
 another attractor and so on. This approach is indeed only a first
 approximation. It is a quite meaningful one because in this
 setting one can explain many phenomena like for instance
 hysteresis cycles, relaxation oscillations, bursting oscillations
 and more complicated alternance of pulsatile and surge patterns
 of coupled neurons populations. It can also serve as a good representation model
 of the so-called ultrastability introduced in Cybernetics by R. Ashby (\cite{Ash})
 (as the
 property that some systems display an adaptation under some external influences
 to switch in a reversible or non-reversible manner from a stable state to another one)
 and used in life sciences by H. Atlan (\cite{At}, \cite{Francoise}).
 \vskip 1pt
 Also, this approach does not take into account canards. Discovered by E. Benoit, Callot, M. and F.
 Diener (\cite{BCDD}) using non-standard analysis, canards are limit
 periodic sets (limit when $\epsilon\to 0$) of the van der Pol
 system:

\begin{equation}
\label{e1}
{\epsilon}\dot{x}=y-f(x)=y-(\frac{x^3}{3}+x^2)
\end{equation}
\begin{equation}
\dot{y}=x-c(\epsilon)
\end{equation}
where $c(\epsilon)$ ranges between some bounds:
\begin{equation}
\label{e2}
c_0+{\rm exp}(-\alpha/\epsilon)<c(\epsilon)<c_0+{\rm exp}(-\beta/\epsilon),
\quad (\alpha>\beta>0).
\end{equation}
The surprising aspect is that a part of the canard coincides with
a repulsive piece of the critical manifold. More recent
contributions of Dumortier and Roussarie (\cite{DR}) yield
new (standard) proofs for the existence of such limit-periodic
sets. There are now several evidences showing the relevance of
this notion to explain experimental facts observed in physiology.
Note that canards occur in two-parameters families in some narrow
set of parameters.\\
The system we present here, although quite simple,
displays a kind of recurrence of the canard effect which
enhances the delay to bifurcating.

\section{Transcritical Dynamical Bifurcation}

The classical transcritical bifurcation occurs when the parameter $\lambda$
in the equation:

\begin{equation}
\label{e3}
\dot{x}=-{\lambda}x+x^2,
\end{equation}
crosses $\lambda=0$. Equation $\ref{e3}$ displays two equilibria, $x=0$
and $x=\lambda$. For $\lambda>0$, $x=0$ is stable and $x=\lambda$
is unstable. After the bifurcation, $\lambda<0$, $x=0$ is stable
and $x=\lambda$ is unstable. The two axis have ``exchanged" their
stability.\\
The terminology ``Dynamical Bifurcation" refers to the situation
where the bifurcation parameter is replaced by a slowly varying
variable. In the case of the transcritical bifurcation, this
yields:

\begin{eqnarray}
\label{e4}
\dot{x}=-yx+x^2\\
\dot{y}=-\epsilon,
\end{eqnarray}
where $\epsilon$ is assumed to be small.\\
This yields
$$\dot{x}=-(-{\epsilon}t+y_0)x+x^2, \quad (y_0=y(0))$$
which is an integrable equation of Bernoulli type. Its solution
displays:
$$x(t)=\frac{x_0{\rm exp}[-Y(t)]}{1-x_0\int_0^t {\rm exp}[-Y(u)]du} \quad (x_0=x(0)$$
$$Y(t)=\int_0^t y(s)ds=\int_0^t (-{\epsilon}s+y_0)ds=-{\epsilon}\frac{t^2}{2}+y_0 t.$$
If we fix an initial data $(x_0,y_0)$, $y_0>0$, $0<x_0<y_0/2$, and we consider
the solution with this initial data we find easily that it takes
time $t=y_0/{\epsilon}$ to reach the axis $y=0$. If $x_0$ is quite
small, that means the orbit stays closer and closer of the
attractive part of the critical manifold untill it reaches the axis $x=y$ and then
coordinate $x$ start increasing. But now consider time
$cy_0/{\epsilon}$, $1\geq c\geq 2$. Then a straightforward computation shows that
$Y(t)=c(1-\frac{c}{2})\frac{y_0^2}{\epsilon}
=\frac{k}{\epsilon},$ and that
$${\rm lim}_{\epsilon\to 0} \int_{\epsilon}^{2y_0/{\epsilon}}{\rm exp}
(\epsilon\frac{t^2}{2}-y_0 t)dt=\frac{2}{y_0}.$$
This shows that if $\epsilon$ is small enough and if the initial
data satifies $x_0<y_0/2$, then:
$$x(t)=O(\frac{x_0{\rm exp}(-k/{\epsilon})}{1-2x_0/y_0})<x_0.$$
This yields that, despite the repulsiveness of the axis $x=0, y<0$,
the orbit remains for a very long time close to $x=0$, indeed
$x(t)<x_0$. Note that after a larger time the orbit blows away
from this repulsive axis. This phenomenon, although quite simply
explained, is of the same nature as the delay to bifurcation
discovered for the dynamical Hopf bifurcation. See for instance
(\cite{ErnBaRinz}, \cite{CDD}, \cite{Ern}, \cite{Z}).
This well-known effect is
instrumental in the example we construct in this article. Some
related work has been done in computing exit points through the
passage of single turning points, for example \cite{Ben},
\cite{DeDu}.

\section{An example of system with enhanced delay}
Consider the equation:
$$\dot{x}=(1-x^2)(x-y)$$
$$\dot{y}={\epsilon}x.$$
The fast dynamics displays the invariant lines $x=-1$, $x=1$, and
$y=x$. A quick analysis shows that, as the slow variable $y$
varies, the fast system undergoes two transcritical bifurcations
near the points $(-1,-1)$ and $(1,1)$. As we recalled in the first
paragraph, a typical orbit near $(-1,-1)$ first displays a ``delay"
along the repulsive part $(x=-1,y<-1)$ of the slow manifold. Then,
by hysteresis, it jumps to the attractive part $(x=1,y<1)$ till it
reaches the other transcritical bifurcation where it again
displays another delay along the repulsive part $(x=1,y>1)$. Then
it jumps again to $(x=-1, y>1)$ and starts again. There is such a
mechanism of successive enhancements of the delay after several
turns generated by the hysteresis. After this intuitive
explanation, we give now a formal proof of the:\\
{\bf Theorem}\\
For all initial data inside the strip $-1<x<1$, for all $\delta$
and for all $T$, the corresponding orbit spends a time larger than
$T$ within a distance less than $\delta$ to the repulsive part of
the slow manifold.\\
{\bf Proof}\\
Inside the strip $\mid x\mid<1$, it is convenient to use the
variable $u$: $x={\rm tanh}u$. The system yields the equations:
$$\dot{u}={\rm tanh}u-y$$
$$\dot{y}={\epsilon}{\rm tanh}u.$$
Consider the function
$$\Phi(u,y)={\frac{1}{2}}({\rm tanh}u-y)^2+{\epsilon}{\rm ln}({\rm cosh}u),$$
and its time derivation along the flow. This displays:
$$\frac{d}{dt}\Phi(u,y)=(\frac{\dot{u}}{{\rm cosh}u})^2$$
Hence the function $\Phi$ is strictly increasing along the flow
(Lyapunov function for the flow).\\
Note now that if $(u(t),y(t))$ is a solution then $(-u(t),-y(t))$
is also a solution. To study the orbits of the system, we can
restrict to initial data $u=u_0\geq 0$ and $y=y_0$.
The first step of the proof is to show that all orbits intersect both axes $u=0$ and
$y=u$ in infinitely many points.\\
Assume first $y_0\geq {\rm tanh}u_0$. Then ${\dot u}(0)\leq 0$. But
$$\frac{d}{dt}(y-{\rm tanh}u)={\epsilon}{\rm tanh}u+\frac{y-{\rm tanh}u}{{\rm cosh}^2u}$$
shows that $y-{\rm tanh}u$ grows hence remains positive. Assume that
$u$ would remain always positive. Then, as $y-{\rm tanh}u>0$, $u$ is
monotone decreasing. Hence there exists $l$ such that $u\to l$ as
$t\to+\infty$.\\
If $l>0$, $\dot{y}={\epsilon}{\rm tanh}u$ implies (via the mean
value theorem) $y\to+\infty$ but then ${\dot u}\to-\infty$ and
(mean value theorem) contradiction with $u>0$.\\
If $l=0$, $\dot{u}+y\to 0$. But $y$ is monotone increasing. If
$y\to+\infty$, then ${\dot u}\to-\infty$ and again contradiction.
If $y$ tends to a finite limit $m$, then $\dot{u}\to -m$ and again
contradiction. Hence all orbits with initial data $(u_0,y_0)$ with
$y_0\geq {\rm tanh}u_0\geq 0$ intersect the axe $u=0$.\\
Consider now the case $y_0<{\rm tanh}u_0$. The variable $u$ is first
strictly increasing (as soon as $y<{\rm tanh}u)$, hence positive and
so $y$ is increasing. Assume that ${\rm tanh}u-y$ would remain
positive along the orbit. Then as $t\to+\infty$, $u$ would tend to
a limit $m$ (eventually $m=+\infty$).As ${\dot y}\to m$, mean
value theorem would imply $y\to +\infty$ and again a contradiction
with ${\dot u}\to-\infty$. So the orbit necessarily intersects the
axe $y={\rm tanh}u$ and ultimately the axe $u=0$ by the preceding
argument. By symmetry, we can also show the existence of two
sequences of times $(t_n)$ and $(\theta_n)$ such that:
$$t_n<\theta_n<t_{n+1}$$
$$x(t_n)=x(t_{n+1})=0$$
$$x(\theta_n)=y(\theta_n)$$
$$y(t_n)=(-1)^na_n,\quad a_n>0.$$
In the second part of the proof we show that the sequence $(a_n)$ is
unbounded.\\
Consider now the function
$$2\Phi(u,y)=w(x,y)=(x-y)^2-{\epsilon}{\rm ln}\mid 1-x^2\mid,$$
which satisfies:
$$\dot{w}=2(1-x^2)(x-y)^2/\epsilon.$$
As $w$ is strictly increasing, this yields:
$$w(t_n)=a_n^2<w(\theta_n)=-{\epsilon}{\rm ln}(1-\xi_n^2)<w(t_{n+1})=a_{n+1}^2,$$
with $\xi_n=x(\theta_n).$
Integration along the flow of $(\dot w)=2\dot{x}(x-y)$ yields:
$$a_{n+1}^2-a_n^2=\int_{t_n}^{t_{n+1}}2\dot{x}(x-y)dt=2\int_{t_n}^{t_{n+1}}x\dot{y}dt$$
$$=2{\epsilon}\int_{t_n}^{t_{n+1}}x^2dt\leq 2{\epsilon}(t_{n}-t_{n+1}).$$
This shows that if $(t_n)$ converges to a finite limit then so does $(a_n)$ and
$(\xi_n)$. Now integration along the flow of $\dot{y}=\epsilon u$
yields:
$$a_n+a_{n+1}\leq {\epsilon}\int_{t_n}^{t_{n+1}}\mid {\rm tanh}u\mid dt\leq
{\epsilon}(t_{n+1}-t_n),$$
and this shows that the sequence of times $(t_n)$ is necessarily
unbounded. Now assume that the sequence $(a_n)$ would be bounded.
Then, as the sequence $(t_n)$ tends to $+\infty$, the function $w$
would be bounded on the orbit. But then so would be both the two
functions $(x-y)^2$ and $-{\rm ln}(1-x^2)$. But then there would
exist a constant $\alpha$ such that $(1-x^2)\geq \alpha$ along the
orbit and
$${\dot w}\geq 2{\alpha}(x-y)^2\geq 2{\alpha}w,$$
hence ${\rm e}^{-2{\alpha}t}w(t)$ increasing and contradiction
with the fact that $w$ would be bounded.\\
Last step is classical in slow-fast dynamics. Consider the
unbounded sequence of points $(0,a_n)$ on the orbit. By Tikhonov's
theorem there is an invariant curve which passes near that
point. By Takens's theorem the flow is conjugated to the fast flow
along this invariant curve. Hence necessarily there is for all
orbits inside the strip, all $\delta$ and all $T$  a part of the
orbit which remains at a distance less than $\delta$ of the
repulsive parts of the boundary of the strip for a time larger
than $T$.

\section{Asymptotics of the system outside the strip}

{\bf Theorem}
\vskip 1pt
Given any initial point $(x_0,y_0)$ outside the strip $\mid x\mid\leq
1$, the corresponding orbit is asymptotic to $y=x$.
\vskip 1pt
{\bf Proof}
\vskip 1pt
We can always assume that $x_0>1$ because the system is symmetric
relatively to the origin. The equations yield:
$$\epsilon\frac{dx}{dy}=(x-\frac{1}{x})(y-x).$$
So if $y\geq x$ and $(x_0,y_0)$, $(x_1,y_1)$ are two points on the
same orbit with $x_0<x_1$, we get:
$$(x_0-\frac{1}{x_0})(y-x)\leq {\epsilon}\frac{dx}{dy}\leq (x_1-\frac{1}{x_1})(y-x).$$
Set:
$$\alpha_i=\frac{1}{\epsilon}(x_i-\frac{1}{x_i}),\quad i=0,1$$
this yields
$$\frac{dx}{dy}-\alpha_0(y-x)\geq 0,$$
$$\frac{dx}{dy}-\alpha_1(y-x)\leq 0,$$
hence
$$\frac{d}{dy}({\rm e}^{\alpha_0y}x)\geq \alpha_0y{\rm e}^{\alpha_0y}$$
$$\frac{d}{dy}({\rm e}^{\alpha_1y}x)\leq \alpha_1y{\rm e}^{\alpha_1y}.$$
Integration between $y_0$ and $y_1$ yields
$${\rm e}^{\alpha_0y_0}(y_0-x_0-\frac{1}{\alpha_0})\geq {\rm e}^{\alpha_0y_1}
(y_1-x_1-\frac{1}{\alpha_0}),$$
and
$${\rm e}^{\alpha_1y_0}(y_0-x_0-\frac{1}{\alpha_1})\leq {\rm e}^{\alpha_1y_1}
(y_1-x_1-\frac{1}{\alpha_1}).$$
The second inequality shows that if $y_0=x_0$, then
$${\rm e}^{\alpha_1y_0}\geq {\rm e}^{\alpha_1y_1}(1-\alpha_1(y_1-x_1)),$$
and thus that the orbit stays above the line $y=x$. If we now
choose $x_0$ (and $\alpha_0)$ large enough, the first inequality
displays:
$$(y_1-x_1)\leq \frac{1}{2\alpha_0}.$$
This shows that the orbit is asymptotic to $y=x$.

\end{document}